\documentclass{amsart}

\usepackage{amsfonts,amssymb,amsmath}

\newcommand{\R}{{\mathbb R}}
\newcommand{\Z}{{\mathbb Z}}
\newcommand{\C}{{\mathbb C}}
\newcommand{\di}{{\rm d}}
\newcommand{\dv}[2]{\langle #1,#2\rangle}

\theoremstyle{plain}
\newtheorem{theorem}{Theorem}[section]
\newtheorem{lemma}[theorem]{Lemma}

\theoremstyle{definition}
\newtheorem{definition}[theorem]{Definition}
\theoremstyle{remark}
\newtheorem{remark}[theorem]{Remark}

\begin{document}

\markboth{A.~G.~Smirnov}{On topological tensor products}

\title{On topological tensor products of functional Fr\'echet and DF spaces}

\author{A.~G.~Smirnov}

\address{I.~E.~Tamm Theory Department, P.~N.~Lebedev Physical Institute,\\ Leninsky prospect 53, Moscow 119991, Russia}
\email{smirnov@lpi.ru}


\begin{abstract}
A convenient technique for calculating completed topological tensor products of functional Fr\'echet and DF spaces is developed. The general
construction is applied to proving kernel theorems for a wide class of weighted spaces of smooth and entire analytic functions.
\end{abstract}

\maketitle

\section{Introduction}
\label{s1}

Let $X$ and $Y$ be sets and $F$, $G$, and $H$ be locally convex spaces consisting of functions defined on $X$, $Y$, and $X\times Y$ respectively. If
the function $(x,y)\to f(x)g(y)$ belongs to $H$ for any $f\in F$ and $g\in G$, then $F\otimes G$ is identified with a linear subspace of $H$. In
applications, it is often important to find out whether $H$ can be interpreted as the completion of $F\otimes G$ with respect to some natural tensor
product topology. Results of this type (or rather their reformulations in terms of bilinear forms) for concrete functional spaces are known as kernel
theorems. For example, Schwartz's kernel theorem for tempered distributions amounts to the statement that the space $S(\R^{k_1+k_2})$ of rapidly
decreasing smooth functions on $\R^{k_1+k_2}$ is identical to the completion $S(\R^{k_1})\widehat{\otimes}_\pi S(\R^{k_2})$ of $S(\R^{k_1})\otimes
S(\R^{k_2})$ with respect to the projective\footnotemark[1] topology. Very general conditions ensuring the \emph{algebraic} coincidence of $H$ with
$F\widehat \otimes_\pi G$ can be derived from the description of completed tensor products of functional spaces given by Grothendieck (\cite{Grot},
Chapitre~2, Th\'eor\`eme~13). Namely, suppose both $F$ and $G$ are Hausdorff and complete, $F$ is nuclear and representable as an inductive limit of
Fr\'echet spaces, and the topology of $F$ is stronger than that of simple convergence. Then $F\widehat \otimes_\pi G$ is algebraically identified
with $H$ if and only if $H$ consists exactly of all functions $h$ on $X\times Y$ such that $h(x,\cdot)\in G$ for every $x\in X$ and the function
$h_v(x)=\dv{v}{h(x,\cdot)}$ belongs to $F$ for all $v\in G'$ (for a locally convex space $G$, we denote by $G'$ and $\dv{\cdot}{\cdot}$ the
continuous dual of $G$ and the canonical bilinear form on $G'\times G$ respectively). In practice, however, the space $H$ usually carries its own
topology and one needs to prove the \emph{topological} coincidence of $H$ and $F\widehat \otimes_\pi G$. In this paper, we find a convenient
criterion (see Theorem~\ref{t1a} below) ensuring such a coincidence under the assumption that $F$ and $G$ are either both Fr\'echet spaces or both
reflexive DF spaces. This criterion allows us to obtain simple proofs of kernel theorems for well-known spaces $K(M_p)$ of smooth functions~\cite{GV}
(or, more precisely, for their suitable generalization, see Definition~\ref{d2s}) and for a class of weighted spaces of entire analytic functions
including those arising in quantum field theory (see, e.g., \cite{Sm, So}).

\section{General construction}
\label{s2}

Given a locally convex space $F$, we denote by $F'_\sigma$, $F'_\tau$, and $F'_b$ the space $F'$ endowed with its weak topology, Mackey topology, and
strong topology respectively. We refer the reader to~\cite{Grot1} for the definition and properties of DF spaces. Here we only mention that the
strong dual of a Fr\'echet space (resp., DF space) is a DF space (resp., Fr\'echet space). Recall that Fr\'echet spaces and their Mackey duals are
B-complete (see~\cite{Schaefer}, Sec.~IV.8). In particular, reflexive DF spaces are B-complete (and hence complete). In Theorem~\ref{t1a} below, we
assume that the considered vector spaces are either real or complex. The ground field ($\R$ or $\C$) is denoted by $\mathbb K$.

\begin{theorem}
 \label{t1a}
Let $X$ and $Y$ be sets, $F$ and $G$ be either both Fr\'echet spaces or both reflexive DF spaces consisting of scalar functions on $X$ and $Y$
respectively, and $H$ be a barelled Hausdorff complete space consisting of scalar functions on $X\times Y$. Let $G$ be nuclear and the topologies of
$F$, $G$, and $H$ be stronger than the topology of simple convergence. Suppose the following conditions are satisfied:
\begin{itemize}
\item[($i$)] For every $f\in F$ and $g\in G$, the function $(x,y)\to f(x)g(y)$ on $X\times Y$ belongs to $H$ and the bilinear mapping $\Phi\colon
F\times G\to H$ taking $(f,g)$ to this function is separately continuous. \item[($ii$)] If $h\in H$, then $h(x,\cdot)\in G$ for every $x\in X$ and
the function $h_v(x)=\dv{v}{h(x,\cdot)}$ belongs to $F$ for every $v\in G'$.
\end{itemize}
Then $\Phi$ is continuous and induces the topological isomorphism $F\widehat{\otimes}_\pi G\simeq H$.
\end{theorem}
\begin{proof}
Recall that separate continuity is equivalent to continuity for bilinear mappings defined on Fr\'echet or barrelled DF spaces (see~\cite{Schaefer},
Theorem~III.5.1 and \cite{Grot1}, Corollaire du Th\'eor\`eme~2). Therefore, $\Phi$ is continuous. Let $\Phi_*\colon F\widehat{\otimes}_\pi G \to H$
be the continuous linear mapping determined by $\Phi$. As usual, let $\mathfrak B_e(F'_\sigma, G'_\sigma)$ denote the space of separately continuous
bilinear forms on $F'_\sigma\times G'_\sigma$ equipped with the biequicontinuous convergence topology (i.e., the topology of the uniform convergence
on the sets of the form $A\times B$, where $A$ and $B$ are equicontinuous sets in $F'$ and $G'$ respectively). Let $S$ be the natural continuous
linear mapping $F\widehat{\otimes}_\pi G \to \mathfrak B_e(F'_\sigma, G'_\sigma)$ which takes $f\otimes g$ to the bilinear form $(u,v)\to
\dv{u}{f}\dv{v}{g}$. Since $G$ is nuclear, $S$ is a topological isomorphism (see~\cite{Grot}, Chapitre~2, Th\'eor\`eme~6 or~\cite{Schaefer}).
Further, let $T$ be the linear mapping\footnotemark[2] $\mathfrak B_e(F'_\sigma, G'_\sigma)\to \mathbb K^{X\times Y}$ defined by the relation
$(Tb)(x,y)=b(\delta_x,\delta_y)$, $b\in \mathfrak B_e(F'_\sigma, G'_\sigma)$ (if $x\in X$ and $y\in Y$, then $\delta_x$ and $\delta_y$ are the linear
functionals on $F$ and $G$ such that $\dv{\delta_x}{f}=f(x)$ and $\dv{\delta_y}{g}=g(y)$; they are continuous because the topologies of $F$ and $G$
are stronger than the topology of simple convergence). Obviously, $T$ is continuous and $TS$ coincides on $F\otimes G$ with $j\Phi_*$, where $j\colon
H\to \mathbb K^{X\times Y}$ is the inclusion mapping. By continuity, we have
\begin{equation}\label{eq1}
j\Phi_*=TS
\end{equation}
everywhere on $F\widehat{\otimes}_\pi G$. Moreover, $T$ is injective because $\delta$-functionals are weakly dense in $F'$ and $G'$. By (\ref{eq1}),
this implies the injectivity of $\Phi_*$.

We now prove the surjectivity of $\Phi_*$. In view of (\ref{eq1}) and the injectivity of $j$ it suffices to show that $\mathrm{Im}\,j\subset
\mathrm{Im}\,j\Phi_*=\mathrm{Im}\,T$. Let $h\in H$ and $L\colon G'_\tau\to F$ be the linear mapping taking $v\in G'$ to $h_v$ (by ($ii$), $L$ is well
defined). We claim that the graph $\mathcal G$ of $L$ is closed. It suffices to show that if an element of the form $(0,f)$ belongs to the closure
$\bar{\mathcal G}$ of $\mathcal G$, then $f=0$. Suppose the contrary that there is $f_0\in F$ such that $f_0\ne 0$ and $(0,f_0)\in \bar{\mathcal G}$.
Let $x_0\in X$ be such that $f_0(x_0)\ne 0$ and let the neighborhood $U$ of $f_0$ be defined by the relation $U=\{f\in F :
|\dv{\delta_{x_0}}{f-f_0}|< |f_0(x_0)|/2\}$. Let $V=\{v\in G' : |\dv{v}{h(x_0,\cdot)}|<|f_0(x_0)|/2\}$. If $f\in U$ and $v\in V$, then we have
$|h_v(x_0)|<|f_0(x_0)|/2<|f(x_0)|$. Hence the neighborhood $V\times U$ of $(0,f_0)$ does not intersect $\mathcal G$. This contradicts the assumption
that $(0,f_0)\in \bar{\mathcal G}$, and our claim is proved. Being nuclear and complete, $G$ is semireflexive and hence $G'_\tau$ is barrelled. We
can therefore apply the closed graph theorem (\cite{Schaefer}, Theorem~IV.8.5) and conclude that $L$ is continuous. Let the bilinear form $b$ on
$F'\times G'$ be defined by the relation $b(u,v)=\dv{u}{h_v}$. The continuity of $L$ implies that $b\in \mathfrak B_e(F'_\sigma, G'_\sigma)$. Since
$b(\delta_x, \delta_y)=h(x,y)$, we have $h=Tb$. Thus, $h\in\mathrm{Im}\,T$ and the surjectivity of $\Phi_*$ is proved.

Thus, $\Phi_*$ is a one-to-one mapping from $F\widehat{\otimes}_\pi G$ onto $H$. If $F$ and $G$ are Fr\'echet (resp., reflexive DF) spaces, then
$F\widehat{\otimes}_\pi G$ is a Fr\'echet (resp., reflexive DF) space (for DF spaces, this follows from \cite{Grot}, Chapitre~2, Th\'eor\`eme~12).
Hence $F\widehat{\otimes}_\pi G$ is B-complete and the open mapping theorem (\cite{Schaefer}, Theorem~IV.8.3) ensures that $\Phi_*$ is a topological
isomorphism.
\end{proof}

\begin{remark}
In contrast to the treatment in~\cite{Treves,Komatsu,BN}, where the density of $F\otimes G$ in $H$ (for concrete functional spaces) was proved ``by
hand'', we have obtained this density automatically as a consequence of nuclearity and the density of $\delta$-functionals in dual spaces. In some
cases (especially for spaces of analytic functions), a direct check of the density of $F\otimes G$ in $H$ may present a considerable difficulty.
\end{remark}

\begin{remark}
Of course, we can exchange the roles of $F$ and $G$ in the formulation of Theorem~\ref{t1a}, i.e., we can assume that $F$ is nuclear and replace
($ii$) by the condition
\begin{itemize}
\item[($ii'$)] If $h\in H$, then $h(\cdot,y)\in F$ for every $y\in Y$ and the function $y\to\dv{u}{h(\cdot,y)}$ belongs to $G$ for every $u\in F'$.
\end{itemize}
Analogous changes can be made in the formulations of Lemmas~\ref{l6s} and~\ref{l6} and Theorems~\ref{t2s} and~\ref{t2} below.
\end{remark}

\section{Kernel theorems for spaces of smooth and analytic functions}
\label{s3}

We now apply Theorem~\ref{t1a} to proving kernel theorems for some spaces of smooth functions. In what follows, we use the standard multi-index
notation:
\[
|\mu|=\mu_1+\ldots+\mu_k,\quad \partial^\mu
f(x)=\frac{\partial^{|\mu|}f(x)}{\partial x_1^{\mu_1}\ldots
\partial x_k^{\mu_k}}\quad (\mu\in \Z_+^k).
\]

\begin{definition}
 \label{d2s}
Let $M=\{M_\gamma\}_{\gamma\in \Gamma}$ be a family of nonnegative measurable functions on $\R^k$ which are bounded on every bounded subset of $\R^k$
and satisfy the following conditions:
\begin{itemize}
\item[(a)] For every $\gamma_1,\gamma_2\in\Gamma$, one can find $\gamma\in\Gamma$ and $C>0$ such that $M_\gamma\geq C (M_{\gamma_1}+M_{\gamma_2})$.
\item[(b)] There is a countable set $\Gamma'\subset \Gamma$ with the property that for every $\gamma\in \Gamma$, one can find $\gamma'\in \Gamma'$
and $C>0$ such that $C M_\gamma\leq M_{\gamma'}$. \item[(c)] For every $x\in \R^k$, one can find $\gamma\in\Gamma$, a neighborhood $O(x)$ of $x$, and
$C>0$ such that $M_\gamma(x')\geq C$ for all $x'\in O(x)$.
\end{itemize}
The space $\mathcal K(M)$ consists of all smooth functions $f$ on $\R^k$ having the finite seminorms
\begin{equation}\label{2s}
\|f\|_{\gamma,m} = \sup_{x\in \R^k,\,|\mu|\leq m}
M_\gamma(x)|\partial^\mu f(x)|
\end{equation}
for all $\gamma\in\Gamma$ and $m\in \Z_+$. The space $\mathcal K_p(M)$, $p\geq 1$, consists of all smooth functions $f$ on $\C^k$ having the finite
seminorms
\begin{equation}\label{3s}
\|f\|^p_{\gamma,m} = \left(\int_{\R^k}
[M_\gamma(x)]^p\,\sum_{|\mu|\leq m}|\partial^\mu f(x)|^p\, \di
x\right)^{1/p}.
\end{equation}
The spaces $\mathcal K(M)$ and $\mathcal K_p(M)$ are endowed with the topologies determined by seminorms~(\ref{2s}) and~(\ref{3s}) respectively.
\end{definition}

We shall say that $M=\{M_\gamma\}_{\gamma\in \Gamma}$ is a defining family of functions on $\R^k$ if it satisfies all requirements of
Definition~\ref{d2s}. Note that if all $M_\gamma$ are strictly positive and continuous, then condition~(c) holds automatically. Condition~(b) ensures
that $\mathcal K(M)$ and $\mathcal K_p(M)$ possess a countable fundamental system of neighborhoods of the origin. It is easy to see that $\mathcal
K(M)$ is actually a Fr\'echet space. Indeed, let $f_n$ be a Cauchy sequence in $\mathcal K(M)$. Then it follows from~(c) that $\partial^\mu f_n(x)$
converge uniformly on every compact subset of $\R^k$ for every multi-index $\mu$. This implies that $f_n$ converge pointwise to a smooth function
$f$. For $\varepsilon>0$, $\gamma\in \Gamma$, and $m\in\Z_+$, choose $n_0$ such that $\|f_{n+l}-f_n\|_{\gamma,m}<\varepsilon$ for all $n\geq n_0$ and
$l\in\Z_+$. Then $M_\gamma(x)|\partial^\mu f_{n+l}(x)-\partial^\mu f_n(x)|<\varepsilon$ for every $x\in\R^k$ and $|\mu|\leq m$. Passing to the limit
$l\to \infty$, we obtain $M_\gamma(x)|\partial^\mu f(x)-\partial^\mu f_n(x)|<\varepsilon$, i.e., $\|f-f_n\|_{\gamma,m}<\varepsilon$ for $n\geq n_0$.
Hence it follows that $f\in \mathcal K(M)$ and $f_n\to f$ in this space.

\begin{lemma}
 \label{ls2}
Let $M=\{M_\gamma\}_{\gamma\in \Gamma}$ be a defining family of functions on $\R^k$. The space $\mathcal D(\R^k)$ of smooth functions with compact
support is dense in $\mathcal K_p(M)$ for any $p\geq 1$.
\end{lemma}

\begin{proof} Let $f\in \mathcal K_p(M)$ and $\varphi\in \mathcal
D(\R^k)$ be such that $\varphi(x)=1$ for $|x|\leq 1$ ($|\cdot|$ is a norm on $\R^k$). For $n=1,2,\ldots$, we define $\varphi_n \in \mathcal D(\R^k)$
by the relation $\varphi_n(x)=\varphi(x/n)$ and set $\psi_n=1-\varphi_n$. To prove the statement, it suffices to show that $\varphi_n f\to f$ (or,
which is the same, that $\psi_n f\to 0$) in $\mathcal K_p(M)$ as $n\to \infty$. Let $\gamma\in\Gamma$, $m\in \Z_+$, and $\mu$ be a multi-index such
that $|\mu|\leq m$. An elementary estimate using the Leibniz formula gives $|\partial^\mu (\psi_n f)(x)|\leq A2^m\sum_{|\nu|\leq m}|\partial^{\nu}
f(x)|$, where $A=1+\sup_{x,\,|\nu|\leq m} |\partial^{\nu}\varphi(x)|$. Since $\psi_n$ vanishes for $|x|\leq n$, it hence follows that
\[
\|\psi_n f\|^p_{\gamma,m}\leq A2^m
q(m)\left(\int_{|x|>n}[M_\gamma(x)]^p\,\sum_{|\mu|\leq
m}|\partial^\mu f(x)|^p\, \di x\right)^{1/p},
\]
where $q(m)$ is the number of multi-indices whose norm does not exceed $m$. Since $\|f\|^p_{\gamma,m}<\infty$, the integral in the right-hand side
tends to zero as $n\to\infty$.
\end{proof}

\begin{lemma}
 \label{l4s}
Let $M=\{M_\gamma\}_{\gamma\in \Gamma}$ be a defining family of functions on $\R^k$ satisfying the following conditions:
\begin{itemize}
\item[$\mathrm{(I)}$] For every $\gamma\in \Gamma$, there are $\gamma'\in \Gamma$ and a bounded integrable nonnegative function $L_{\gamma\gamma'}$
on $\R^k$ such that $M_\gamma\leq L_{\gamma\gamma'}M_{\gamma'}$ and $L_{\gamma\gamma'}(x)\to 0$ as $|x|\to\infty$.

\item[$\mathrm{(II)}$] For every $\gamma\in \Gamma$, there are $\gamma'\in \Gamma$, a neighborhood of the origin $B$ in $\R^k$, and $C>0$ such that
$M_\gamma(x)\leq CM_{\gamma'}(x+y)$ for any $x\in \R^k$ and $y\in B$.
\end{itemize}
Then the space $\mathcal K(M)$ is nuclear and coincides, both as a set and topologically, with $\mathcal K_p(M)$ for all $p\geq 1$.
\end{lemma}

\begin{proof} Let $f\in \mathcal K(M)$, $\gamma\in \Gamma$, and
$m\in \Z_+$. Choosing $\gamma'$ and $L_{\gamma\gamma'}$ such that (I) is satisfied, we obtain $\|f\|^p_{\gamma,m}\leq A\|f\|_{\gamma',m}$, where
$A=\left(q(m)\int [L_{\gamma\gamma'}(x)]^p\,\di x\right)^{1/p}<\infty$ (as above, $q(m)$ is the number of multi-indices with the norm $\leq m$).
Hence we have a continuous inclusion $\mathcal K(M)\subset\mathcal K_p(M)$.

We now prove that the topology induced on $\mathcal K(M)$ from $\mathcal K_p(M)$ coincides with the original topology of $\mathcal K(M)$. In other
words, given $\gamma\in\Gamma$ and $m\in \Z_+$, we have to find $\tilde\gamma\in\Gamma$, $\tilde m\in\Z_+$, and $A>0$ such that
\begin{equation}\label{7s}
\|f\|_{\gamma,m}\leq A\|f\|^p_{\tilde\gamma,\tilde m}
\end{equation}
for every $f\in \mathcal K(M)$. By (II), there are $\gamma',\gamma''\in \Gamma$, a neighborhood of the origin $B\subset \R^k$, and $C>0$ such that
$M_\gamma(x)\leq CM_{\gamma'}(x+y)$ and $M_{\gamma'}(x)\leq CM_{\gamma''}(x+y)$ for any $x\in \R^k$ and $y\in B$. Let $\psi$ be a smooth nonnegative
function such that $\int \psi(x)\di x=1$ and $\mathrm{supp}\,\psi\subset B$. Set $\tilde M(x)=\int M_{\gamma'}(x+x') \psi(x') \di x'$. Then $\tilde
M$ is a smooth function on $\R^k$ and we have
\begin{align}
& M_\gamma(x)=\int M_\gamma(x)\psi(x')\di x'\leq C\int
M_{\gamma'}(x+x')\psi(x')\di x'=C\tilde M(x),\label{4s}\\
& |\partial^\mu \tilde M(x)|=\left|\int
M_{\gamma'}(x+x')\partial^\mu\psi(x')\di x'\right|\leq C_\mu
M_{\gamma''}(x),\label{5s}
\end{align}
where $C_\mu = C\int |\partial^\mu\psi(x)|\di x$. In view of condition~(I) inequality~(\ref{5s}) implies that $|\partial^{\nu}\tilde M(x)\partial^\mu
f(x)|\to 0$ as $|x|\to\infty$ for all multi-indices $\mu$ and $\nu$ and every $f\in \mathcal K(M)$. For $|\mu|\leq m$, it hence follows
from~(\ref{4s}) and~(\ref{5s}) that
\begin{multline}
M_\gamma(x)|\partial^\mu f(x)| \leq C\tilde M(x) |\partial^\mu
f(x)|=\\= C\left|\int_{-\infty}^{x_1}\di
x'_1\ldots\int_{-\infty}^{x_k}\di x'_k\frac{\partial^k}{\partial
x'_1\ldots\partial x'_k} [\tilde M(x') \partial^\mu
f(x')]\right|\leq C'\|f\|^1_{\gamma'',m+k},\label{6s}
\end{multline}
where $C'=C\sum_{|\mu|\leq k} C_\mu$. Let $\tilde m= m+k$ and $\tilde \gamma$ be such that $M_{\gamma''}\leq L_{\gamma'',\tilde\gamma}\,M_{\tilde
\gamma}$, where $L_{\gamma'',\,\tilde\gamma}(x)$ is integrable and tends to zero as $|x|\to\infty$. Estimating $\|f\|^1_{\gamma'',m+k}$ by the
H\"older inequality, we conclude from~(\ref{6s}) that (\ref{7s}) holds with
\[
A=C'\left(q(\tilde
m)\int[L_{\gamma'',\tilde\gamma}(x)]^{p/(p-1)}\,\di
x\right)^{(p-1)/p}.
\]
Since $\mathcal D(\R^k)\subset \mathcal K(M)$, it follows from Lemma~\ref{ls2} that $\mathcal K(M)$ is a dense subspace of $\mathcal K_p(M)$. At the
same time, the completeness of $\mathcal K(M)$ implies that it is closed in $\mathcal K_p(M)$. Hence, we have $\mathcal K(M)=\mathcal K_p(M)$.

To prove the nuclearity of $\mathcal K(M)$, we shall use the following criterion obtained by Pietsch~\cite{Pietsch}.

\begin{lemma}
 \label{l5}
A locally convex space $F$ is nuclear if and only if some (every) fundamental system $\mathcal U$ of absolutely convex neighborhoods of the origin
has the following property:

For every neighborhood of the origin $U\in\mathcal U$, there is a neighborhood of the origin $V\in\mathcal U$ and a positive Radon
measure\footnotemark[3] $\tau$ on $V^\circ$ such that
\begin{equation}\label{4}
p_U(f)\leq\int_{V^\circ}|\dv{u}{f}|\,\di\tau(u)
\end{equation}
for every $f\in F$ ($p_U$ is the Minkowski functional of the set $U$).
\end{lemma}

For $\gamma\in \Gamma$ and $m\in\Z_+$, we set $U_{\gamma,m}=\{f\in \mathcal K(M) : \|f\|_{\gamma,m}\leq 1\}$. By condition~(a) of
Definition~\ref{d2s}, the scalar multiples of $U_{\gamma,m}$ form a fundamental system of neighborhoods of the origin in $\mathcal K(M)$, and we have
$p_{U_{\gamma,m}}(f)=\|f\|_{\gamma,m}$. Fix $\gamma$ and $m$ and choose $\tilde \gamma$, $\tilde m$, and $A$ such that inequality~(\ref{7s}) with
$p=1$ is satisfied for all $f\in \mathcal K(M)$. Let $\tilde \gamma'$ be such that $M_{\tilde\gamma}\leq
L_{\tilde\gamma,\tilde\gamma'}M_{\tilde\gamma'}$, where $L_{\tilde\gamma,\tilde\gamma'}$ is integrable. Acting as in the derivation of
formulas~(\ref{4s}) and~(\ref{5s}), we find a smooth function $\tilde M$, an index $\tilde\gamma''$, and $C>0$ such that $M_{\tilde\gamma'}\leq
C\tilde M$ and $\tilde M\leq C M_{\tilde\gamma''}$. For every $x\in\R^k$ and multi-index $\mu$, we define the functional $\varepsilon^\mu_x\in
\mathcal K'(M)$ by the relation
\[
\dv{\varepsilon^\mu_x}{f}=\tilde M(x)
\partial^\mu f(x)/C,\quad f\in \mathcal K(M).
\]
If $|\mu|\leq\tilde m$, then we obviously have $|\dv{\varepsilon^\mu_x}{f}|\leq 1$ for $\|f\|_{\tilde\gamma'',\tilde m}\leq 1$, i.e.,
$\varepsilon^\mu_x\in U^\circ_{\tilde\gamma'',\tilde m}$. Moreover, the mapping $x\to\varepsilon^\mu_x$ from $\R^k$ to $\mathcal K'(M)$ is weakly
continuous. Hence the function $\varphi(\varepsilon^\mu_x)$ is bounded and continuous on $\R^k$ for every continuous function $\varphi$ on the weakly
compact set $U^\circ_{\tilde\gamma'',\tilde m}$. Therefore, the formula
\[
\tau(\varphi) = A C^2\int_{\R^k} L_{\tilde\gamma,\tilde\gamma'}(x)
\sum_{|\mu|\leq\tilde m}\varphi(\varepsilon^\mu_x)\,\di x,\quad
\varphi\in C(U^\circ_{\tilde\gamma'',\tilde m}),
\]
defines a positive Radon measure $\tau$ on $U^\circ_{\tilde\gamma'',\tilde m}$. It follows from this definition that
\[
\|f\|_{\gamma,m}\leq A C^2\int_{\R^k}
L_{\tilde\gamma,\tilde\gamma'}(x) \sum_{|\mu|\leq \tilde
m}|\dv{\varepsilon^\mu_x}{f}| \,\di x=
\int_{U^\circ_{\tilde\gamma'',\tilde m}} |\dv{u}{f}|\,\di\tau(u)
\]
for every $f\in \mathcal K(M)$. In view of Lemma~\ref{l5} this estimate shows that $\mathcal K(M)$ is nuclear. Lemma~\ref{l4s} is proved.
\end{proof}

\noindent\emph{Examples.} (1) Let $\Gamma$ be the set of all compact subsets of $\R^k$ and $M_\gamma$ be the characteristic function of $\gamma$.
Then $\mathcal K(M)$ is the space $C^\infty(\R^k)$ endowed with its standard topology. Conditions~(I) and~(II) are obviously satisfied.

\medskip
\noindent (2) Let $\Gamma=\Z_+$ and $M_l = (1+|x|)^l$. Then $\mathcal K(M)$ is the Schwartz space $S(\R^k)$ of rapidly decreasing functions.
Conditions~(I) and~(II) are obviously satisfied.

\medskip
\noindent (3) Let $\alpha>0$, $A\geq 0$, $\Gamma$ be the interval $(A,\infty)$, and $M_{A'}(x)=\exp(|x/A'|^{1/\alpha})$ for every $A'>A$. Then
$\mathcal K(M)$ coincides with the Gelfand--Shilov space $S_{\alpha,\tilde A}$, where $\tilde A=(\alpha/e)^\alpha A$ (see~\cite{GS}, Section~IV.3).
Conditions~(I) and~(II) are obviously satisfied and, therefore, the space $S_{\alpha,A}$ is nuclear for any $A\geq 0$.

\begin{remark} \label{rrrr}
The spaces $\mathcal K(M)$ are similar to the spaces $K(M_p)$ introduced in the  classical book~\cite{GV}. Theorem~I.7 in~\cite{GV} asserts that
condition~(I) of Lemma~$\ref{l4s}$ (which is called condition~(N) there) is sufficient for the nuclearity of $K(M_p)$. However, in our opinion, the
proof of this theorem is not quite satisfactory (an estimate of type~(\ref{7s}) is obtained there with a constant $A$ depending implicitly on the
function $f$). Moreover, the condition $M_p(x)\geq 1$ included in the definition of $K(M_p)$ is actually redundant.
\end{remark}

Let $M=\{M_\gamma\}_{\gamma\in \Gamma}$ and $N=\{N_\omega\}_{\omega\in \Omega}$ be defining families of functions on $\R^{k_1}$ and $\R^{k_2}$
respectively. We denote by $M\otimes N$ the family formed by the functions
\[
(M\otimes N)_{\gamma\omega}(x,y)=M_\gamma(x) N_\omega(y),\quad
(\gamma,\omega)\in \Gamma\times \Omega.
\]
Clearly, $M\otimes N$ is a defining family of functions on $\R^{k_1+k_2}$.

\begin{lemma}
 \label{l6s}
Let $M=\{M_\gamma\}_{\gamma\in \Gamma}$ and $N=\{N_\omega\}_{\omega\in \Omega}$ be defining families of functions on $\R^{k_1}$ and $\R^{k_2}$
respectively and let $h\in \mathcal K(M\otimes N)$. Suppose $N$ satisfies conditions {\rm (I)} and {\rm (II)} of Lemma~$\ref{l4s}$. Then
$h(x,\cdot)\in \mathcal K(N)$ for every $x\in \R^{k_1}$ and the function $h_v(x)=\dv{v}{h(x,\cdot)}$ belongs to $\mathcal K(M)$ for all $v\in
\mathcal K'(N)$. Moreover, for every multi-index $\mu\in\Z_+^{k_1}$, we have
\begin{equation}\label{diff}
\partial^\mu h_v(x)= \dv{v}{\partial^\mu_x h(x,\cdot)}.
\end{equation}
\end{lemma}

\begin{proof} By Lemma~\ref{l4s}, $\mathcal K(N)$ is a nuclear
Fr\'echet space. This implies, in particular, that it is reflexive. Let $Q(M)$ be the space consisting of the sequences $\psi=\{\psi^\mu\}_{\mu\in
\Z_+^{k_1}}$ of functions on $\R^{k_1}$ having the finite norms
\[
|||\psi|||_{\gamma,m}=\sup_{x\in\R^{k_1},\,|\mu|\leq m}
|\psi^\mu(x)| M_\gamma(x)
\]
for all $\gamma\in\Gamma$ and $m\in \Z_+$. Let $T\colon \mathcal K(M)\to Q(M)$ be the mapping taking $f\in \mathcal K(M)$ to the sequence
$\{\partial^\mu f\}$. Obviously, $T$ maps $\mathcal K(M)$ isomorphically onto its image, and since $\mathcal K(M)$ is complete, $\mathrm{Im}\,T$ is a
closed subspace of $Q(M)$. For $v\in \mathcal K'(N)$ and $\mu\in\Z_+^{k_1}$, we set $\psi^\mu_v(x)= \dv{v}{\partial^\mu_x h(x,\cdot)}$. Since
$h_v=\psi_v^\mu$ for zero $\mu$, it suffices to show that the sequence $\psi_v=\{\psi_v^\mu\}$ belongs to $\mathrm{Im}\,T$. For every
$\omega\in\Omega$ and $n\in \Z_+$, we set $B_{\omega,n}=\{v\in \mathcal K'(N): |\dv{v}{g}|\leq \|g\|_{\omega,n} \,\,\forall g\in \mathcal K(N)\}$. If
$\gamma\in \Gamma$, $m\in \Z_+$, $|\mu|\leq m$, and $v\in B_{\omega,n}$, then we have
\[
|\dv{v}{\partial^\mu_x h(x,\cdot)}| M_\gamma(x)\leq
\|\partial^\mu_x h(x,\cdot)\|_{\omega,n}M_\gamma(x)\leq
\|h\|_{(\gamma,\omega),\,m+n},\quad x\in \R^{k_1}.
\]
Hence, $|||\psi_v|||_{\gamma,m}\leq \|h\|_{(\gamma,\omega),\,m+n}$ for $v\in B_{\omega,n}$. Thus, $\psi_v$ belongs to the space $Q(M)$ for any $v\in
\mathcal K'(N)$ and the image of $B_{\omega,n}$ under the mapping $v\to \psi_v$ is bounded in $Q(M)$. The scalar multiples of $B_{\omega,n}$ form a
fundamental system of bounded subsets in the space $\mathcal K_b'(N)$, which is bornological as the strong dual of a reflexive Fr\'echet space
(\cite{Schaefer}, Corollary~1 of Proposition~IV.6.6). Therefore, the mapping $v\to \psi_v$ from $\mathcal K_b'(N)$ to $Q(M)$ is continuous. If
$v=\delta_y$ for some $y\in \R^{k_2}$, then $\psi_v$ obviously belongs to $\mathrm{Im}\,T$. This implies that $\psi_v\in \mathrm{Im}\,T$ for all
$v\in\mathcal K'(N)$ because $\mathrm{Im}\,T$ is closed in $Q(M)$, the linear span of $\delta$-functionals is dense in $\mathcal K_b'(N)$ by the
reflexivity of $\mathcal K(N)$, and the image of the closure of a set under a continuous mapping is contained in the closure of the image of this
set.
\end{proof}

We are now ready to prove the kernel theorem for the spaces $\mathcal K(M)$.
\begin{theorem}
 \label{t2s}
Let $M$ and $N$ be defining families of functions on $\R^{k_1}$ and $\R^{k_2}$ respectively and the bilinear mapping $\Phi\colon \mathcal K(M)\times
\mathcal K(N)\to \mathcal K(M\otimes N)$ be defined by the relation $\Phi(f,g)(x,y)=f(x)g(y)$. If $N$ satisfies conditions {\rm (I)} and {\rm (II)}
of Lemma~$\ref{l4s}$, then $\Phi$ induces the topological isomorphism $\mathcal K(M)\widehat{\otimes}_\pi \mathcal K(N)\simeq \mathcal K(M\otimes
N)$.
\end{theorem}

\begin{proof}
Lemma~\ref{l4s} implies that $\mathcal K(N)$ is a nuclear Fr\'echet space. The required statement therefore follows from Theorem~\ref{t1a} (the
fulfilment of $(i)$ is obvious and $(ii)$ is ensured by Lemma~\ref{l6s}).
\end{proof}

We now consider the spaces of entire analytic functions. We say that a family $M=\{M_\gamma\}_{\gamma\in \Gamma}$ of functions on $\C^k$ is a
defining family of functions on $\C^k$ if $M$ is a defining family of functions on the underlying real space $\R^{2k}$. In what follows, we identify
defining families of functions on $\C^k$ with the corresponding defining families of functions on $\R^{2k}$. In particular, if $M$ is a defining
family of functions on $\C^k$, then $\mathcal K(M)$ will denote the corresponding space of $C^\infty$-functions on $\R^{2k}$ and the statement that
$M$ satisfies conditions~(I) and~(II) of Lemma~$\ref{l4s}$ will mean that (I) and~(II) are fulfilled if $M$ is viewed as a family of functions on
$\R^{2k}$. If $M$ and $N$ are defining families of functions on $\C^{k_1}$ and $\C^{k_2}$ respectively, then $M\otimes N$ will be interpreted as a
defining family of functions on $\C^{k_1+k_2}$.

\begin{definition}
 \label{d2}
Let $M=\{M_\gamma\}_{\gamma\in \Gamma}$ be a defining family of functions on $\C^k$. The space $\mathcal H(M)$ consists of all entire analytic
functions $f$ on $\C^k$ having the finite seminorms
\begin{equation}\label{2}
\|f\|_\gamma = \sup_{z\in \C^k} M_\gamma(z)|f(z)|.
\end{equation}
For $p\geq 1$, the space $\mathcal H_p(M)$ consists of all entire analytic functions $f$ on $\C^k$ having the finite seminorms
\begin{equation}\label{3}
\|f\|^p_\gamma = \left(\int_{\C^k} [M_\gamma(z)]^p\,|f(z)|^p\, \di
\lambda(z)\right)^{1/p},
\end{equation}
where $\di\lambda$ is the Lebesgue measure on $\C^k$. The spaces $\mathcal H(M)$ and $\mathcal H_p(M)$ are endowed with the topologies determined by
the seminorms~(\ref{2}) and~(\ref{3}) respectively.
\end{definition}

The same arguments as in the case of $\mathcal K(M)$ show that $\mathcal H(M)$ is a Fr\'echet space for any defining family of functions $M$.

\begin{lemma}
 \label{l4}
Let $M=\{M_\gamma\}_{\gamma\in \Gamma}$ be a defining family of functions on $\C^k$ satisfying conditions~{\rm (I)} and~{\rm (II)} of
Lemma~$\ref{l4s}$. Then the space $\mathcal H(M)$ is nuclear and coincides, both as a set and topologically, with $\mathcal H_p(M)$ for all $p\geq
1$.
\end{lemma}

\begin{proof}
Let $\tilde{\mathcal H}(M)$ be the subspace of $\mathcal K(M)$ consisting of all elements of $\mathcal K(M)$ which are entire analytic functions.
Since a subspace of a nuclear space is nuclear, it follows from Lemma~\ref{l4s} that $\tilde{\mathcal H}(M)$ is a nuclear space. Therefore, to prove
the nuclearity of $\mathcal H(M)$, it suffices to show that $\tilde {\mathcal H}(M)=\mathcal H(M)$. We obviously have the continuous inclusion
$\tilde {\mathcal H}(M)\subset\mathcal H(M)$. Conversely, let $f\in \mathcal H(M)$, $\gamma\in \Gamma$, and $\mu,\nu\in \Z_+^k$ be multi-indices. By
(II), we can find $\gamma'\in \Gamma$, a neighborhood of the origin $B\subset \C^k$, and $C>0$ such that $M_\gamma(z)\leq CM_{\gamma'}(z+z')$ for any
$z\in \C^k$ and $z'\in B$. For $r>0$ and $z\in \C^k$, let $D_r(z)$ denote the polydisk with the radius $r$ centered at $z$. If $\zeta-z\in B$, then
we have $M_\gamma(z)|f(\zeta)|\leq C\|f\|_{\gamma'}$. Therefore, choosing $r>0$ so small that $D_r(0)\subset B$ and using the Cauchy formula, we
obtain
\begin{multline}
|\partial^{\mu}_x\partial^{\nu}_y
f(x+iy)|M_\gamma(x+iy)=|\partial^{\mu+\nu}_z f(z)|M_\gamma(z)\leq
\\ \leq \frac{(\mu+\nu)!}{(2\pi)^k}\oint\limits_{|\zeta_1-z_1|=r} \di\zeta_1\ldots
\oint\limits_{|\zeta_k-z_k|=r}\di\zeta_k \left|\frac{M_\gamma(z)
f(\zeta)}{(\zeta-z)^{\iota+\mu+\nu}}\right|\leq C (\mu+\nu)!\,
r^{-|\mu+\nu|}\|f\|_{\gamma'},\nonumber
\end{multline}
where $z=x+iy$ and $\iota$ is the multi-index $(1,\ldots,1)$. This inequality implies that $\|f\|_{\gamma,m}\leq C m!\, r^{-m}\|f\|_{\gamma'}$ for
any $m\in\Z_+$ ($\|f\|_{\gamma,m}$ is given by~(\ref{2s})). Hence $\mathcal H(M)$ is continuously embedded in $\tilde{\mathcal H}(M)$ and, therefore,
$\tilde{\mathcal H}(M)=\mathcal H(M)$.

We now prove the coincidence of $\mathcal H(M)$ and $\mathcal H_p(M)$. Let $f\in \mathcal H(M)$ and $\gamma\in \Gamma$. By~(I), we can find
$\gamma'\in \Gamma$ such that $M_\gamma\leq L_{\gamma\gamma'}M_{\gamma'}$, where $L_{\gamma\gamma'}(z)$ is integrable and tends to zero as
$|z|\to\infty$. Then we obtain $\|f\|^p_\gamma\leq A\|f\|_{\gamma'}$, where $A=\left(\int [L_{\gamma\gamma'}(z)]^p\,\di
\lambda(z)\right)^{1/p}<\infty$. Hence we have the continuous inclusion $\mathcal H(M)\subset\mathcal H_p(M)$. Conversely, let $f\in \mathcal H_p(M)$
and $\gamma\in \Gamma$. Let $\gamma'$, $B$, $C$, and $r$ be as in the preceding paragraph. It follows from the Cauchy formula that $f(z)=(\pi
r^2)^{-k}\int_{D_r(z)} f(\zeta)\di\lambda(\zeta)$ for every $z\in \C^k$. Multiplying both parts of this relation by $M_\gamma(z)$ and using~(II) and
the H\"older inequality, we obtain
\[
|f(z)|M_\gamma(z) \leq C\left(\frac{1}{(\pi
r^2)^k}\int_{D_r(z)}[M_{\gamma'}(\zeta)]^p\,|f(\zeta)|^p\,
\di\lambda(\zeta) \right)^{1/p}.
\]
Extending integration in the right-hand side to the whole of $\C^k$ and passing to the supremum in the left-hand side, we find that $\|f\|_\gamma\leq
C(\pi r^2)^{-k/p}\|f\|^p_{\gamma'}$. Hence $\mathcal H_p(M)$ is continuously embedded in $\mathcal H(M)$ and, therefore, $\mathcal H(M)=\mathcal
H_p(M)$.
\end{proof}

\begin{lemma}
 \label{l6}
Let $M=\{M_\gamma\}_{\gamma\in \Gamma}$ and $N=\{N_\omega\}_{\omega\in \Omega}$ be defining families of functions on $\C^{k_1}$ and $\C^{k_2}$
respectively and let $h\in \mathcal H(M\otimes N)$. Suppose $N$ satisfies conditions~{\rm (I)} and~{\rm (II)} of Lemma~$\ref{l4s}$. Then
$h(z,\cdot)\in \mathcal H(N)$ for every $z\in \C^{k_1}$ and the function $h_v(z)=\dv{v}{h(z,\cdot)}$ belongs to $\mathcal H(M)$ for all $v\in
\mathcal H'(N)$.
\end{lemma}

\begin{proof}
Let $v\in \mathcal H'(N)$. As shown in the proof of Lemma~\ref{l4}, $\mathcal H(N)$ is the subspace of $\mathcal K(N)$ consisting of those elements
of $\mathcal K(N)$ that are entire analytic functions. By the Hahn--Banach theorem, $v$ has a continuous extension $\hat v$ to $\mathcal K(N)$. Then
$h_v(z)=\dv{\hat v}{h(z,\cdot)}$ and Lemma~\ref{l6s} implies that $h_v\in \mathcal K(M)$. Moreover, it follows from~(\ref{diff}) that $h_v$ satisfies
the Cauchy--Riemann equations and, therefore, is an entire analytic function. Hence $h_v\in \mathcal H(M)$ and the lemma is proved.
\end{proof}

We are now ready to prove the kernel theorem for the spaces $\mathcal H(M)$.
\begin{theorem}
 \label{t2}
Let $M$ and $N$ be defining families of functions on $\C^{k_1}$ and $\C^{k_2}$ respectively and the bilinear mapping $\Phi\colon \mathcal H(M)\times
\mathcal H(N)\to \mathcal H(M\otimes N)$ be defined by the relation $\Phi(f,g)(x,y)=f(x)g(y)$. If $N$ satisfies conditions {\rm (I)} and {\rm (II)}
of Lemma~$\ref{l4s}$, then $\Phi$ induces the topological isomorphism $\mathcal H(M){\widehat{\otimes}}_\pi \mathcal H(N)\simeq \mathcal H(M\otimes
N)$.
\end{theorem}

\begin{proof} Lemma~\ref{l4} implies that $\mathcal H(N)$
is a nuclear Fr\'echet space. The required statement therefore follows from Theorem~\ref{t1a} (the fulfilment of $(i)$ is obvious and $(ii)$ is
ensured by Lemma~\ref{l6}).
\end{proof}

\begin{remark}
Instead of using Theorem~~\ref{t1a} and Lemma~\ref{l6}, one can obtain Theorem~\ref{t2} by combining the results of~\cite{Bier} with Lemma~\ref{l4}.
More precisely, it follows from Proposition~3.5 in~\cite{Bier} that $\mathcal H(M){\widehat{\otimes}}_\varepsilon \mathcal H(N)\simeq \mathcal
H(M\otimes N)$, where $\mathcal H(M){\widehat{\otimes}}_\varepsilon \mathcal H(N)$ is the completion of $\mathcal H(M)\otimes \mathcal H(N)$ with
respect to the biequicontinuous topology. In view of the nuclearity of $\mathcal H(N)$ ensured by Lemma~\ref{l4} this implies Theorem~\ref{t2}.
\end{remark}

{\small

\par\medskip\noindent
{\bf Acknowledgements.} The author is grateful to M.A.~Soloviev for useful discussions. The research was partially supported by the grants
RFBR~05-02-17451 and LSS-1615.2008.2.

\par\bigskip\noindent
{\bf Notes.}
\begin{itemize}
\item[1.] Recall that the projective (inductive) topology on $F\otimes G$ is the strongest locally convex topology on $F\otimes G$ such that the
canonical bilinear mapping $(f,g)\to f\otimes g$ is continuous (resp., separately continuous). In general, one must carefully distinguish between the
inductive and projective topologies. However, these topologies always coincide for Fr\'echet and barrelled DF spaces which are considered in this
paper. For definiteness, we speak everywhere about the projective topology.

\item[2.] The space $\mathbb K^{X\times Y}$ of all scalar functions on $X\times Y$ is endowed with the topology of simple convergence or, which is
the same, the direct product topology.

\item[3.] Recall that a Radon measure on a compact set $K$ is, by definition, a continuous linear form on the space $C(K)$ of continuous functions on
$K$. Recall also that the polar set $V^\circ$ of a neighborhood of the origin $V$ in a locally convex space is weakly compact.
\end{itemize}
}

\end{document}